\documentclass[a4paper,10pt,fleqn]{article}
\usepackage[latin1]{inputenc}
\usepackage[left=2cm,right=2cm,bottom=3cm,top=3cm,]{geometry}

\usepackage[all,2cell,v2,arrow, curve, frame, matrix, graph, tips]{xy}

\UseAllTwocells

\usepackage[square,numbers,sort]{natbib}
\usepackage{tikz-cd}
\usepackage{ifthen}
\usepackage{bbm} 
\usepackage{amsfonts}
\usepackage{savesym}
\usepackage{mathtools}
\usepackage{amsmath}
\usepackage{amssymb}
\usepackage{tcolorbox}
\usepackage{enumerate}
\usepackage[thmmarks,standard]{ntheorem}
\usepackage{microtype}
\usepackage{stmaryrd}

\usepackage{verbatim}
\usepackage[bookmarks=true,
            bookmarksnumbered=true, breaklinks=true,
            pdfstartview=FitH, hyperfigures=false,
            plainpages=false, naturalnames=true,
            colorlinks=true,pagebackref=true,
            pdfpagelabels]{hyperref}
\usepackage{cleveref}
\usepackage[toc,page]{appendix}
\theoremstyle{plain}

\newcommand{\E}{\ensuremath{\mathbb{S}\text{ets}}}

\newcommand{\I}{\ensuremath{\infty\text{-}}}

\newcommand{\Mc}{\ensuremath{\infty\text{-}\mathbb{M}\mathbb{C}\text{at}}}

\newcommand{\Mmc}{\ensuremath{(\infty,m)\text{-}\mathbb{M}\mathbb{C}\text{at}}}

\newcommand*{\G}{\ensuremath{\mathbb{G}\text{r}}}

\newcommand*{\TG}{\ensuremath{T\text{-}\G}}

\newcommand*{\MEt}{\ensuremath{\infty\text{-}\mathbb{M}\mathbb{E}\text{tC}}}

\newcommand*{\MEtm}{\ensuremath{(\infty,m)\text{-}\mathbb{M}\mathbb{E}\text{tC}}}
\newcommand{\N}{\ensuremath{\mathbb{N}}}

\title{Notes on Multiple Higher Category Theory} \author{Camell
  Kachour}
\begin{document}
\maketitle
\vspace*{3.5cm}
\begin{abstract}
These notes follows the articles \cite{kamel, Cam, cam-cubique} which show how powerful can be the method of
\textit{Stretchings} initiated with the \textit{Globular Geometry} by Jacques Penon in \cite{penon} ,
to weakened \textit{strict higher structures}. Here we adapt this method to weakened strict
multiple $\infty$-categories, strict multiple $(\infty,m)$-categories, and in particular we obtain algebraic models
of weak multiple $\infty$-groupoids.

\end{abstract}

\begin{minipage}{118mm}{\small
      }
\end{minipage}

\hypersetup{%
  linkcolor=blue}%
\tableofcontents
\vspace*{1cm}

\vspace*{1cm}

\section*{Introduction}
\textit{Strict multiple categories} had been introduced by
Charles Ehresman in \cite{Ehresmann} in order to produce higher generalization of categories. Surprisingly the \textit{multiple geometry} used in it to produce a theory of higher category has not been studied as much as it deserve. In these notes we hope filling this gap where in particular we use the technology of \textit{Stretchings}
to produce algebraic models of multiple higher category theory. More
specifically we shall introduce :
\begin{itemize}

\item Algebraic models of weak multiple $\infty$-categories in \ref{weak}

\item Algebraic models of weak multiple $(\infty,m)$-categories in \ref{weakm}
\end{itemize}

In particular we propose algebraic models of weak multiple
$\infty$-groupoids which are our models of weak multiple $(\infty,0)$-categories.


\section{Multiple Sets}
\subsection{Multiple Sets}
Fix an integer $n\geq 1$, each finite sequence $\underline{i}=(i_1,i_2,...,i_n)\in\N^n$ such that
$1\leq i_1<i_2<...<i_n$ is called an $n$-color, and all $n$-colors form a set denoted $\mathfrak{I}_n$.
Also for a fix $n$-color $\underline{i}\in\mathfrak{I}_n$, $S^1_{\underline{i}}=\{i_1,...,i_n\}$ denotes its underlying set
of $1$-colors
$i_j$. If $1\leq k\leq n$, the set $S^k_{\underline{i}}$ of $k$-colors of $\underline{i}$ is well understood, it has
$C^n_k=\frac{n!}{k!(n-k)!}$ $k$-colors :
$S^k_{\underline{i}}=\{\underline{i}^k_{1},\underline{i}^k_{2},...,\underline{i}^k_{C^n_k}\}$, i.e each $\underline{i}^k_{j}$
is a subsequence $(i_{j_1},i_{j_1},...,i_{j_k})\in\N^k$ of $\underline{i}$ such that
$1\leq i_{j_1}<i_{j_2}<...<i_{j_k}$. The set $S^{n-1}_{\underline{i}}$ of $(n-1)$-colors of $\underline{i}$ has
a particular importance : for $S^{n-1}_{\underline{i}}=\{\underline{i}^{n-1}_{1},\underline{i}^{n-1}_{2},...,\underline{i}^{n-1}_{n}\}$
and $1\leq k\leq n$ we use the notation $\underline{i}^{n-1}_{k}=\underline{i}-i_k$
which means that the $(n-1)$-color $\underline{i}^{n-1}_{k}$
is $(i_1,i_2,...,\hat{i_k},...,i_n)\in\N^{n-1}$ i.e we delete the
$1$-color $i_k$ from $\underline{i}=(i_1,i_2,...,i_n)\in\N^n$. In fact if $\underline{i}=(i_1,i_2,...,i_n)\in\N^n$
is an $n$-color and if $l$ is an integer such that the $(n+1)$-sequence $(i_1,i_2,...,l,...,i_n)$ is an
$(n+1)$-color, then we write this last $(n+1)$-color by $\underline{i}+l$. Thus if $\underline{i}=(i_1,i_2,...,i_n)\in\N^n$
is an $n$-color then the notation "\textit{minus}" $\underline{i}-i_k$ means the $(n-1)$-color $(i_1,i_2,...,\hat{i_k},...,i_n)\in\N^{n-1}$
and the notation "\textit{add}" $\underline{i}+l$ means that $l$ doesn't belongs to $S^1_{\underline{i}}$ and that $\underline{i}+l$
must be seen as its corresponding $(n+1)$-color, and possibly it can be reindexed if necessary.

An $n$-colored set means a
set $C_{\underline{i}}$ colored by an $n$-color $\underline{i}$ like just above.
An $n$-multiple data $C_n$ means a countable set $C_n$ of $n$-colored sets
$C_n=(C_{\underline{i}})_{\underline{i}\in\mathfrak{I}_n}$, and a multiple data
$C$ means a countable set $C=(C_n)_{n\in\N}$ of $n$-multiple datas $C_n$
if $n\geq 1$ and a \textit{set of objects} $C_0$.

A multiple set is given by a multiple data $C=(C_n)_{n\in\N}$ such that for all
$n\geq 1$, and all $n$-colored set $C_{\underline{i}}$ in it, this $n$-colored set
is equipped for all $i_j\in S^1_{\underline{i}}$ with \textit{sources and targets} :

$$\begin{tikzcd}
C_{\underline{i}}\arrow[rr,"s^{\underline{i}}_{\underline{i}-i_j}"]&&C_{\underline{i}-i_j}
\end{tikzcd},\qquad\begin{tikzcd}
C_{\underline{i}}\arrow[rr,"t^{\underline{i}}_{\underline{i}-i_j}"]&&C_{\underline{i}-i_j}
\end{tikzcd}$$

such that for any $n$-color $\underline{i}$ with $n\geq 2$ the following diagrams commute :

$$
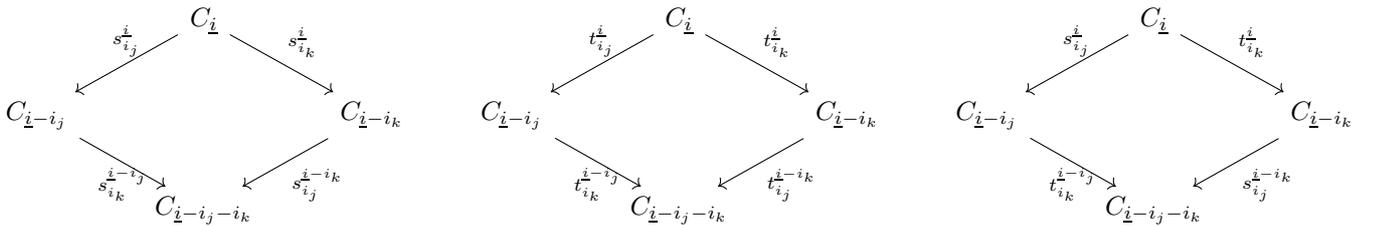
\begin{tikzcd}
&C_{\underline{i}}\arrow[ld,"s^{\underline{i}}_{i_j}"{above}]\arrow[rd,"s^{\underline{i}}_{i_k}"]\\
C_{\underline{i}-i_j}\arrow[rd,"s^{\underline{i}-i_j}_{i_k}"{below}]&&C_{\underline{i}-i_k}\arrow[ld,"s^{\underline{i}-i_k}_{i_j}"]\\
&C_{\underline{i}-i_j-i_k}
\end{tikzcd}\qquad\begin{tikzcd}
&C_{\underline{i}}\arrow[ld,"t^{\underline{i}}_{i_j}"{above}]\arrow[rd,"t^{\underline{i}}_{i_k}"]\\
C_{\underline{i}-i_j}\arrow[rd,"t^{\underline{i}-i_j}_{i_k}"{below}]&&C_{\underline{i}-i_k}\arrow[ld,"t^{\underline{i}-i_k}_{i_j}"]\\
&C_{\underline{i}-i_j-i_k}
\end{tikzcd}\qquad\begin{tikzcd}
&C_{\underline{i}}\arrow[ld,"s^{\underline{i}}_{i_j}"{above}]\arrow[rd,"t^{\underline{i}}_{i_k}"]\\
C_{\underline{i}-i_j}\arrow[rd,"t^{\underline{i}-i_j}_{i_k}"{below}]&&C_{\underline{i}-i_k}\arrow[ld,"s^{\underline{i}-i_k}_{i_j}"]\\
&C_{\underline{i}-i_j-i_k}
\end{tikzcd}$$

We shall often use the short notation $(C,s,t)$ to denote a multiple set $C$ with sources $s$ and
targets $t$, and no confusion should appear with the langage of colors. If $(C',s',t')$ is
another multiple set, a morphism of multiple sets :

$$\begin{tikzcd}
(C,s,t)\arrow[rrr,"f"]&&&(C',s',t')
\end{tikzcd}$$
is given for all $n$-color $\underline{i}$ by a map $f_{\underline{i}}$ of $\E$

$$\begin{tikzcd}
C_{\underline{i}}\arrow[rrr,"f_{\underline{i}}"]&&&C'_{\underline{i}}
\end{tikzcd}$$

which is compatible with all sources and all targets :

$$\begin{tikzcd}
C_{\underline{i}}\arrow[dd,"s^{\underline{i}}_{i_j}"{left}]
\arrow[rrr,"f_{\underline{i}}"]&&&C'_{\underline{i}}\arrow[dd,"s'^{\underline{i}}_{i_j}"]\\\\
C_{\underline{i}-i_j}\arrow[rrr,"f_{\underline{i}-i_j}"{below}]&&&C'_{\underline{i}-i_j}
\end{tikzcd}\qquad \begin{tikzcd}
C_{\underline{i}}\arrow[dd,"t^{\underline{i}}_{i_j}"{left}]
\arrow[rrr,"f_{\underline{i}}"]&&&C'_{\underline{i}}\arrow[dd,"t'^{\underline{i}}_{i_j}"]\\\\
C_{\underline{i}-i_j}\arrow[rrr,"f_{\underline{i}-i_j}"{below}]&&&C'_{\underline{i}-i_j}
\end{tikzcd}$$
The category of multiple sets is denoted $\mathbb{M}\E$

\subsection{Reflexive Multiple Sets}

A reflexive multiple set $(C,s,t,1)$ is given by a multiple set $(C,s,t)$ such that for all $n$-color $\underline{i}=(i_1,...,i_n)\in\mathfrak{I}_n$ ($n\geq 1$)
and all $1$-color $i_j\in S^1_{\underline{i}}$ the $n$-colored set $C_{\underline{i}}$ is equipped with \textit{reflexions} $1^{\underline{i}-i_j}_{\underline{i}}$

$$\begin{tikzcd}
C_{\underline{i}-i_j}\arrow[rrrr,"1^{\underline{i}-i_j}_{\underline{i}}"]&&&&C_{\underline{i}}
\end{tikzcd}$$

For $n\geq 2$ consider an $n$-color $\underline{i}=(i_1,...,i_n)$ and $j, k\in\llbracket 1,...,n\rrbracket$ such that $j\ne k$, then we require the
following commutative diagrams :

$$\begin{tikzcd}
&C_{\underline{i}-i_j}\arrow[ld,"s^{\underline{i}-i_j}_{i_k}"{above}]\arrow[rd,"1^{\underline{i}-i_j}_{\underline{i}}"]\\
C_{\underline{i}-i_j-i_k}\arrow[rd,"1^{\underline{i}-i_j-i_k}_{\underline{i}-i_k}"{below}]&&C_{\underline{i}}\arrow[ld,"s^{\underline{i}}_{i_k}"]\\
&C_{\underline{i}-i_k}
\end{tikzcd}\qquad\begin{tikzcd}
&C_{\underline{i}-i_j}\arrow[ld,"t^{\underline{i}-i_j}_{i_k}"{above}]\arrow[rd,"1^{\underline{i}-i_j}_{\underline{i}}"]\\
C_{\underline{i}-i_j-i_k}\arrow[rd,"1^{\underline{i}-i_j-i_k}_{\underline{i}-i_k}"{below}]&&C_{\underline{i}}\arrow[ld,"t^{\underline{i}}_{i_k}"]\\
&C_{\underline{i}-i_k}
\end{tikzcd}\qquad\begin{tikzcd}
&C_{\underline{i}}\\
C_{\underline{i}-i_j}\arrow[ru,"1^{\underline{i}-i_j}_{\underline{i}}"]&&C_{\underline{i}-i_k}\arrow[lu,"1^{\underline{i}-i_k}_{\underline{i}}"{above}]\\
&C_{\underline{i}-i_j-i_k}\arrow[lu,"1^{\underline{i}-i_j-i_k}_{\underline{i}-i_j}"]\arrow[ru,"1^{\underline{i}-i_j-i_k}_{\underline{i}-i_k}"{right}]
\end{tikzcd}$$
If $(C',s',t',1')$ is another reflexive multiple set, then a morphism of reflexive multiple
sets :
$$\begin{tikzcd}
(C,s,t;1)\arrow[rrr,"f"]&&&(C',s',t';1')
\end{tikzcd}$$
is given by a morphism $f$ in $\mathbb{M}\E$ which is compatible the reflexivity's structures, that
is for all $n$-color $\underline{i}$ and all $1$-color $i_j\in S^1_{\underline{i}}$, the following
diagram commutes :
$$\begin{tikzcd}
C_{\underline{i}}\arrow[rrr,"f_{\underline{i}}"]&&&C'_{\underline{i}}\\\\
C_{\underline{i}-i_j}\arrow[uu,"1^{\underline{i}-i_j}_{\underline{i}}"]
\arrow[rrr,"f_{\underline{i}-i_j}"{below}]&&&C'_{\underline{i}-i_j}\arrow[uu,"1^{\underline{i}-i_j}_{\underline{i}}"{right}]
\end{tikzcd}$$
The category of reflexive multiple sets is denoted $\mathbb{M}_r\E$

The first important monad of this article is the monad of reflexive
multiple sets given by the forgetful functor

$$\begin{tikzcd}
\mathbb{M}_r\E\arrow[rrr,"U"]&&&\mathbb{M}\E
\end{tikzcd}$$

and which is in fact monadic. The proof that $U$ is right adjoint and
monadic come from the underlying projective sketches of $\mathbb{M}_r\E$ and $\mathbb{M}\E$, and by the evident applications of the theorem of Foltz in
\cite{foltz} and Lair in \cite{lair}.

\section{Strict multiple $\infty$-categories}
\label{strict}

A multiple $\infty$-magma
$(M,s,t,(\circ^{\underline{i}}_{\underline{i}-i_j})_{n\geq 1, \underline{i}\in\mathfrak{I}_n)}$ is
given by a multiple set $(M,s,t)$ such that for all $n\geq 1$ and all $\underline{i}\in\mathfrak{I}_n$
its underlying $n$-colored sets $M_{\underline{i}}$ are equipped with operations
$$\begin{tikzcd}
M_{\underline{i}}\underset{M_{\underline{i}-i_j}}\times M_{\underline{i}}\arrow[rrr,"\circ^{\underline{i}}_{i_j}"]&&& M_{\underline{i}}
\end{tikzcd}$$
where $M_{\underline{i}}\underset{M_{\underline{i}-i_j}}\times M_{\underline{i}}$ are given by the
following pullbacks
$$\begin{tikzcd}
C_{\underline{i}}\underset{C_{\underline{i}-i_j}}\times C_{\underline{i}} \arrow[ddd,"\pi_0"{left}]\arrow[rrr,"\pi_1"]&&&C_{\underline{i}}\arrow[ddd,"t^{\underline{i}}_{i_j}"]\\\\\\
C_{\underline{i}}\arrow[rrr,"s^{\underline{i}}_{i_j}"{below}]&&&C_{\underline{i}-i_j}
\end{tikzcd}$$
and such that these operations $\circ^{\underline{i}}_{i_j}$ follow the following
\textit{positional axioms}

\begin{itemize}
\item $s^{\underline{i}}_{\underline{i}-i_j}(a\circ^{\underline{i}}_{i_j} b)=s^{\underline{i}}_{\underline{i}-i_j}(b)$
 and $t^{\underline{i}}_{\underline{i}-i_j}(a\circ^{\underline{i}}_{i_j} b)=t^{\underline{i}}_{\underline{i}-i_j}(a)$

 \item $s^{\underline{i}}_{\underline{i}-i_k}(a\circ^{\underline{i}}_{i_j} b)=s^{\underline{i}}_{\underline{i}-i_k}(a)
\circ^{\underline{i}}_{i_j}s^{\underline{i}}_{\underline{i}-i_k}(b)$ and
$t^{\underline{i}}_{\underline{i}-i_k}(a\circ^{\underline{i}}_{i_j} b)=t^{\underline{i}}_{\underline{i}-i_k}(a)
\circ^{\underline{i}}_{i_j}t^{\underline{i}}_{\underline{i}-i_k}(b)$ if $j\ne k$
\end{itemize}
A multiple $\infty$-magma shall be denoted with the shorter notation $(M,s,t;\circ)$ when
no confusion appears. If $(M',s',t';\circ)$ is another multiple $\infty$-magma, a morphism
of multiple $\infty$-magmas
$$\begin{tikzcd}
(M,s,t;\circ)\arrow[rrr,"f"]&&&(M',s',t';\circ')
\end{tikzcd}$$
is given by a morphism $f\in\mathbb{M}\E$ such which respects the operations $\circ$, that
is for all $n$-color $\underline{i}$ and all $1$-color $i_j\in S^1_{\underline{i}}$ we have
$f_{\underline{i}}(a\circ^{\underline{i}}_{i_j}b)=f_{\underline{i}}(a)\circ^{\underline{i}}_{i_j}f_{\underline{i}}(b)$. The category of multiple $\infty$-magmas is denoted $\infty\text{-}\mathbb{M}\mathbb{M}\text{ag}$.

A reflexive multiple $\infty$-magma $(M,s,t,1,(\circ^{\underline{i}}_{\underline{i}-i_j})_{n\geq 1, \underline{i}\in\mathfrak{I}_n)}$ is given by a reflexive multiple set $(M,s,t,1)$ and a multiple
$\infty$-magma $(M,s,t,(\circ^{\underline{i}}_{\underline{i}-i_j})_{n\geq 1, \underline{i}\in\mathfrak{I}_n)}$
such that
$$1^{\underline{i}}_{\underline{i}+i_k}(a\circ^{\underline{i}}_{i_j} b)
=1^{\underline{i}}_{\underline{i}+i_k}(a)\circ^{\underline{i}+i_k}_{i_j}1^{\underline{i}}_{\underline{i}+i_k}(b)$$

Morphisms of reflexive multiple $\infty$-magmas are those of $\infty\text{-}\mathbb{M}\mathbb{M}\text{ag}$
which are also morphisms of $\mathbb{M}_r\E$. The category of reflexive multiple $\infty$-magmas
is denoted $\infty\text{-}\mathbb{M}_r\mathbb{M}\text{ag}$.

A strict multiple $\infty$-categories is given by a reflexive multiple $\infty$-magma
$(C,s,t;(\circ^{\underline{i}}_{\underline{i}-i_j})_{n\geq 1, \underline{i}\in\mathfrak{I}_n)}$
such that operations $\circ^{\underline{i}}_{\underline{i}-i_j}$ are associative, are unital
i.e if for all
$a\in C_{\underline{i}}$ we have

$$a\circ^{\underline{i}}_{\underline{i}-i_j}1^{\underline{i}-i_j}_{\underline{i}}(t^{\underline{i}}_{\underline{i}-i_j}(a))
=1^{\underline{i}-i_j}_{\underline{i}}(s^{\underline{i}}_{\underline{i}-i_j}(a))\circ^{\underline{i}}_{\underline{i}-i_j}
=a$$
and follow the \textit{middle-four interchange axiom}
$$(a\circ^{\underline{i}}_{\underline{i}-i_j}b)\circ^{\underline{i}}_{\underline{i}-i_k}
(c\circ^{\underline{i}}_{\underline{i}-i_j}d)=(a\circ^{\underline{i}}_{\underline{i}-i_k}c)\circ^{\underline{i}}_{\underline{i}-i_j}
(b\circ^{\underline{i}}_{\underline{i}-i_k}d)$$
when compositions in both side are well defined.

Morphisms of strict multiple $\infty$-categories are those of $\infty\text{-}\mathbb{M}_r\mathbb{M}\text{ag}$. The category of strict multiple $\infty$-categories is denoted $\Mc$.

The second important monad of this article is the monad $\mathbb{S}=(S,\lambda,\nu)$ of strict multiple $\infty$-categories given by the forgetful functor

$$\begin{tikzcd}
\Mc\arrow[rrr,"U"]&&&\mathbb{M}\E
\end{tikzcd}$$

and which is in fact monadic. The proof that $U$ is right adjoint and
monadic come from the underlying projective sketches of $\Mc$ and $\mathbb{M}\E$ which are not difficult to be described. A similar
sketch is described in \cite{cam-cubique}. Then these results come from an application of the theorem of Foltz in \cite{foltz} and Lair in \cite{lair}.

\section{Weak multiple $\infty$-categories}
\label{weak}

\subsection{Multiple categorical stretchings}
\label{stretchings}

A multiple categorical stretching
$\mathbb{E}=(M,C,
(\pi_{\underline{i}})
_{n\geq 1,\underline{i}\in\mathfrak{I}_n},
([-,-]^{\underline{i}}_{\underline{i}+i_r})
_{n\geq 1,\underline{i}\in\mathfrak{I}_n})$
is given by the following datas :

\begin{itemize}

\item A reflexive multiple $\infty$-magma
$(M,s,t,1,(\circ^{\underline{i}}_{\underline{i}-i_j})_{n\geq 1, \underline{i}\in\mathfrak{I}_n)}$

\item A strict multiple $\infty$-category
$(C,s,t,1,(\circ^{\underline{i}}_{\underline{i}-i_j})_{n\geq 1, \underline{i}\in\mathfrak{I}_n)}$

\item A morphism of reflexive multiple $\infty$-magmas
\begin{tikzcd}
M\arrow[rrrr,"\pi"]&&&&C
\end{tikzcd}

\item Operations $[-,-]^{\underline{i}}_{\underline{i}+i_r}$ :

$$\begin{tikzcd}
\underline{M}_{\underline{i}} \arrow[rrrr,"{[-;-]}^{\underline{i}}_{\underline{i}+i_r}"]&&&&M_{\underline{i}+i_r}
\end{tikzcd}$$
where $\underline{M}_{\underline{i}}=\{(\alpha,\beta)\in M_{\underline{i}}\times M_{\underline{i}} : \pi_{\underline{i}}(\alpha)=\pi_{\underline{i}}(\beta) \}$ and such that :

\begin{itemize}

\item $s^{\underline{i}+i_r}_{i_s}[\alpha;\beta]^{\underline{i}}_{\underline{i}+i_r}=[s^{\underline{i}}_{i_s}(\alpha),s^{\underline{i}}_{i_s}(\beta)]^{\underline{i}-i_s}_{\underline{i}-i_s+i_r}$
and $t^{\underline{i}+i_r}_{i_s}[\alpha;\beta]^{\underline{i}}_{\underline{i}+i_r}=[t^{\underline{i}}_{i_s}(\alpha),t^{\underline{i}}_{i_s}(\beta)]^{\underline{i}-i_s}_{\underline{i}-i_s+i_r}$

\item $s^{\underline{i}+i_j}_{i_j}([\alpha;\beta]^{\underline{i}}_{\underline{i}+i_j})=\alpha$ and $t^{\underline{i}+i_j}_{i_j}([\alpha;\beta]^{\underline{i}}_{\underline{i}+i_j})=\beta$

\item $\pi_{\underline{i}+i_r}([\alpha;\beta]^{\underline{i}}_{\underline{i}+i_r})
=1^{\underline{i}}_{\underline{i}+i_r}(\pi_{\underline{i}}(\alpha))=1^{\underline{i}}_{\underline{i}+i_r}(\pi_{\underline{i}}(\beta))$

\end{itemize}

\end{itemize}

A morphism of multiple categorical stretchings
$$\begin{tikzcd}
\mathbb{E}\arrow[rrr,"{(m,c)}"]&&&\mathbb{E}'
\end{tikzcd}$$

is given by the following commutative square in
$\I\mathbb{M}_{\text{r}}\mathbb{M}\text{ag}$,

$$\begin{tikzcd}
M\arrow[rrr,"m"]
\arrow[dd,"\pi'"{left}]
&&&M'\arrow[dd,"\pi'"]\\\\
C \arrow[rrr,"c"{below}]&&&C'
\end{tikzcd}$$

thus we also have the following square in $\E$ for all $n$-color $\underline{i}\in\mathfrak{I}_n$
$(n\geq 1)$

$$\begin{tikzcd}
M_{\underline{i}}\arrow[rrr,"m_{\underline{i}}"]
\arrow[dd,"\pi_{\underline{i}}'"{left}]
&&&M'_{\underline{i}}\arrow[dd,"\pi'_{\underline{i}}"]\\\\
C_{\underline{i}} \arrow[rrr,"c_{\underline{i}}"{below}]&&&C'_{\underline{i}}
\end{tikzcd}$$

and we require for all $(\alpha,\beta)\in\underline{M}_{\underline{i}}$ the following
equality
$$m_{\underline{i}+i_r}([\alpha,\beta]^{\underline{i}}_{\underline{i}+i_r})=
[m_{\underline{i}}(\alpha),m_{\underline{i}}(\beta)]'^{\underline{i}}_{\underline{i}+i_r}$$

The category of multiple categorical stretchings is denoted $\MEt$.

Now consider the forgetful functor :

$$\begin{tikzcd}
\MEt\arrow[rrr,"U"]&&&\mathbb{M}\E
\end{tikzcd}$$

given by :

$$\begin{tikzcd}
(M,C,
(\pi_{\underline{i}})
_{n\geq 1,\underline{i}\in\mathfrak{I}_n},
([-,-]^{\underline{i}}_{\underline{i}+i_r})
_{n\geq 1,\underline{i}\in\mathfrak{I}_n})\arrow[rrr,mapsto]&&&M
\end{tikzcd}$$

\begin{proposition}
 The functor $U$ just above
 has a left adjoint which produces a monad
$\mathbb{M}=(M,\eta,\nu)$ on the category of multiple sets.
\end{proposition}
The proof that $U$ is right adjoint comes from the underlying projective sketches of $\MEt$ and $\mathbb{M}\E$ which are not difficult to be described. A similar
sketch is described in \cite{cam-cubique}. Then these results come from an application of the theorem of Foltz in \cite{foltz}.

\begin{definition}
A weak multiple $\infty$-category is a $\mathbb{M}$-algebra
\end{definition}

%
%
%

\section{Multiple $(\infty,m)$-Sets}
\label{reverses}

Consider a multiple set $(C,s,t)$, and an $n$-color $\underline{i}$ and a
$1$-color $i_j\in S^1_{\underline{i}}$. A $(\underline{i},i_j)$-reversor on it
is given by a map

$$\xymatrix{C_{\underline{i}}\ar[rrr]^{j^{\underline{i}}_{i_j}}&&&C_{\underline{i}}}$$

such that the following two diagrams commute :

$$\xymatrix{C_{\underline{i}}\ar[rr]^{j^{\underline{i}}_{i_j}}
\ar[rd]_{s^{\underline{i}}_{i_j}}&&C_{\underline{i}}
\ar[ld]^{t^{\underline{i}}_{i_j}}\\
&C_{\underline{i}-i_j}}\qquad
\xymatrix{C_{\underline{i}}\ar[rr]^{j^{\underline{i}}_{i_j}}
\ar[rd]_{t^{\underline{i}}_{i_j}}&&C_{\underline{i}}
\ar[ld]^{s^{\underline{i}}_{i_j}}\\
&C_{\underline{i}-i_j}}$$

If for all $k>m$, and all $k$-color $\underline{i}$, and each $i_j\in S^1_{\underline{i}}$
they are such $(\underline{i},i_j)$-reversor $j^{\underline{i}}_{i_j}$ on $(C,s,t)$ then
we say that it is a multiple $(\infty,m)$-set. The family of maps
$(j^{\underline{i}}_{i_j})_{k>m, \underline{i}\in\mathfrak{I}_k, i_j\in S^1_{\underline{i}}}$
is called a multiple $(\infty,m)$-structure and in that case we shall say that $(C,s,t)$ is equipped with the multiple $(\infty,m)$-structure $j=(j^{\underline{i}}_{i_j})_{k>m, \underline{i}\in\mathfrak{I}_k, i_j\in S^1_{\underline{i}}}$. Seen as multiple $(\infty,m)$-set we denote it by
$(C,s,t;(j^{\underline{i}}_{i_j})_{k>m, \underline{i}\in\mathfrak{I}_k, i_j\in S^1_{\underline{i}}})$,
or just $(C,s,t;j)$ for a shorter notation, where $j$ design its underlying multiple $(\infty,m)$-structure.

If $(C',s',t';(j'^{\underline{i}}_{i_j})_{k>m, \underline{i}\in\mathfrak{I}_k, i_j\in S^1_{\underline{i}}})$
is another multiple $(\infty,m)$-set, then a morphism of multiple $(\infty,m)$-sets
\begin{tikzcd}
(C,s,t;(j^{\underline{i}}_{i_j})_{k>m, \underline{i}\in\mathfrak{I}_k, i_j\in S^1_{\underline{i}}})
\arrow[rrr,"f"]&&&(C',s',t';(j'^{\underline{i}}_{i_j})_{k>m, \underline{i}\in\mathfrak{I}_k, i_j\in S^1_{\underline{i}}})
\end{tikzcd}
is given by a morphism of multiple sets such that for each $k>m$, each $k$-color
$\underline{i}$ and each $i_j\in S^1_{\underline{i}}$ we have the following commutative diagrams
$$\begin{tikzcd}
C_{\underline{i}}
\arrow[rrr,"j^{\underline{i}}_{i_j}"]
\arrow[dd,"f_{\underline{i}}"{left}]
&&&C_{\underline{i}}\arrow[dd,"f_{\underline{i}}"]\\\\
C'_{\underline{i}}\arrow[rrr,"j^{\underline{i}}_{i_j}"{below}]&&&C'_{\underline{i}}
\end{tikzcd}$$

The category of multiple $(\infty,m)$-sets is denoted $(\infty,m)$-$\mathbb{M}\E$

\begin{remark}
\label{main-remark}
The $(\infty,m)$-structures that we used to define multiple $(\infty,m)$-sets have globular
and cubical analogues (see \cite{Cam,cam-cubique}) that we called the \textit{minimal $(\infty,m)$-structures}. The
multiple analogue of the \textit{globular maximal $(\infty,m)$-structures} as defined in \cite{Cam}
and of \textit{cubical maximal $(\infty,m)$-structures} as defined in \cite{cam-cubique}
is as follow : for all $k>m$, for each $k$-color $\underline{i}$ and for all $(k-m-1)$-color
$\{i_{j_1},...,i_{j_{k-m-1}}\}\in S^{k-m-1}_{\underline{i}}$, there exist maps

$$\begin{tikzcd}
C_{\underline{i}}\arrow[rr,"j^{\underline{i}}_{i_{j_1}}"]&&C_{\underline{i}}
\end{tikzcd}, \qquad\begin{tikzcd}
C_{\underline{i}-i_{j_1}}\arrow[rr,"j^{\underline{i}-i_{j_1}}_{i_{j_2}}"]&&C_{\underline{i}-i_{j_1}}
\end{tikzcd}, \qquad\begin{tikzcd}
C_{\underline{i}-i_{j_1}-i_{j_2}}\arrow[rr,"j^{\underline{i}-i_{j_1}-i_{j_2}}_{i_{j_3}}"]&&
C_{\underline{i}-i_{j_1}-i_{j_2}}
\end{tikzcd},$$

$$\begin{tikzcd}
C_{\underline{i}-i_{j_1}-...-i_{j_{k-m-2}}}\arrow[rrr,"j^{\underline{i}-i_{j_1}-...-i_{j_{k-m-2}}}_{i_{j_{k-m-1}}}"]
&&&
C_{\underline{i}-i_{j_1}-...-i_{j_{k-m-2}}}
\end{tikzcd}, \qquad\begin{tikzcd}
C_{\underline{i}-i_{j_1}-...-i_{j_{k-m-1}}}\arrow[rrr,"j^{\underline{i}-i_{j_1}-...-i_{j_{k-m-1}}}_{i_{j_{k-m}}}"]
&&&
C_{\underline{i}-i_{j_1}-...-i_{j_{k-m-1}}}
\end{tikzcd}$$

such that we have the following diagrams in $\E$ which commute serially :

    $$
  \xymatrixrowsep{8mm}
  \xymatrixcolsep{4mm}
  \xymatrix{
  C_{\underline{i}}\ar[rrrr]^{j^{\underline{i}}_{i_{j_1}}}
  \ar[d]<+1pt>_{s^{\underline{i}}_{i_{j_1}}}
  \ar[d]<+5pt>^{t^{\underline{i}}_{i_{j_1}}}&&&&C_{\underline{i}}
  \ar[d]<+3pt>^{s^{\underline{i}}_{i_{j_1}}}
  \ar[d]<+6pt>_{t^{\underline{i}}_{i_{j_1}}}\\
                                C_{\underline{i}-i_{j_1}}\ar[rrrr]^{j^{\underline{i}-i_{j_1}}_{i_{j_2}}}
    \ar[d]<+1pt>_{s^{\underline{i}-i_{j_1}}_{i_{j_2}}}
    \ar[d]<+5pt>^{t^{\underline{i}-i_{j_1}}_{i_{j_2}}}&&&&
      C_{\underline{i}-i_{j_1}}
      \ar[d]<+3pt>^{s^{\underline{i}-i_{j_1}}_{i_{j_2}}}
      \ar[d]<+6pt>_{t^{\underline{i}-i_{j_1}}_{i_{j_2}}}\\
                                C_{\underline{i}-i_{j_1}-i_{j_2}}
                                \ar[rrrr]^{j^{\underline{i}-i_{j_1}-i_{j_2}}_{i_{j_3}}}\ar@{.>}[d]<+1pt>^{}\ar@{.>}[d]<+5pt>^{}&&&&C_{\underline{i}-i_{j_1}-i_{j_2}}\ar@{.>}[d]<+3pt>^{}\ar@{.>}[d]<+6pt>_{}\\
                                C_{\underline{i}-i_{j_1}-...-i_{j_{k-m-2}}}
                                \ar[rrrr]^{j^{\underline{i}-i_{j_1}-...-i_{j_{k-m-2}}}_{i_{j_{k-m-1}}}}
  \ar[d]<+1pt>_{s^{\underline{i}-i_{j_1}-...-i_{j_{k-m-2}}}_{i_{j_{k-m-1}}}}
  \ar[d]<+5pt>^(.4){t^{\underline{i}-i_{j_1}-...-i_{j_{k-m-2}}}_{i_{j_{k-m-1}}}}&&&&
                                C_{\underline{i}-i_{j_1}-...-i_{j_{k-m-2}}}
  \ar[d]<+3pt>^{s^{\underline{i}-i_{j_1}-...-i_{j_{k-m-2}}}_{i_{j_{k-m-1}}}}
  \ar[d]<+6pt>_{t^{\underline{i}-i_{j_1}-...-i_{j_{k-m-2}}}_{i_{j_{k-m-1}}}}\\
                                C_{\underline{i}-i_{j_1}-...-i_{j_{k-m-1}}}
                                \ar[rrrr]^{j^{\underline{i}-i_{j_1}-...-i_{j_{k-m-1}}}_{i_{j_{k-m}}}}
                                \ar[drr]_{t^{\underline{i}-i_{j_1}-...-i_{j_{k-m-1}}}_{i_{j_{k-m}}}}&&&&
                                C_{\underline{i}-i_{j_1}-...-i_{j_{k-m-1}}}
                                \ar[dll]^{s^{\underline{i}-i_{j_1}-...-i_{j_{k-m-1}}}_{i_{j_{k-m}}}}\\
                                &&C_{\underline{i}-i_{j_1}-...-i_{j_{k-m}}} } $$

$$
  \xymatrixrowsep{8mm}
  \xymatrixcolsep{4mm}
  \xymatrix{
  C_{\underline{i}}\ar[rrrr]^{j^{\underline{i}}_{i_{j_1}}}
  \ar[d]<+1pt>_{s^{\underline{i}}_{i_{j_1}}}
  \ar[d]<+5pt>^{t^{\underline{i}}_{i_{j_1}}}&&&&C_{\underline{i}}
  \ar[d]<+3pt>^{s^{\underline{i}}_{i_{j_1}}}
  \ar[d]<+6pt>_{t^{\underline{i}}_{i_{j_1}}}\\
                                C_{\underline{i}-i_{j_1}}\ar[rrrr]^{j^{\underline{i}-i_{j_1}}_{i_{j_2}}}
    \ar[d]<+1pt>_{s^{\underline{i}-i_{j_1}}_{i_{j_2}}}
    \ar[d]<+5pt>^{t^{\underline{i}-i_{j_1}}_{i_{j_2}}}&&&&
      C_{\underline{i}-i_{j_1}}
      \ar[d]<+3pt>^{s^{\underline{i}-i_{j_1}}_{i_{j_2}}}
      \ar[d]<+6pt>_{t^{\underline{i}-i_{j_1}}_{i_{j_2}}}\\
                                C_{\underline{i}-i_{j_1}-i_{j_2}}
                                \ar[rrrr]^{j^{\underline{i}-i_{j_1}-i_{j_2}}_{i_{j_3}}}\ar@{.>}[d]<+1pt>^{}\ar@{.>}[d]<+5pt>^{}&&&&C_{\underline{i}-i_{j_1}-i_{j_2}}\ar@{.>}[d]<+3pt>^{}\ar@{.>}[d]<+6pt>_{}\\
                                C_{\underline{i}-i_{j_1}-...-i_{j_{k-m-2}}}
                                \ar[rrrr]^{j^{\underline{i}-i_{j_1}-...-i_{j_{k-m-2}}}_{i_{j_{k-m-1}}}}
  \ar[d]<+1pt>_{s^{\underline{i}-i_{j_1}-...-i_{j_{k-m-2}}}_{i_{j_{k-m-1}}}}
  \ar[d]<+5pt>^(.4){t^{\underline{i}-i_{j_1}-...-i_{j_{k-m-2}}}_{i_{j_{k-m-1}}}}&&&&
                                C_{\underline{i}-i_{j_1}-...-i_{j_{k-m-2}}}
  \ar[d]<+3pt>^{s^{\underline{i}-i_{j_1}-...-i_{j_{k-m-2}}}_{i_{j_{k-m-1}}}}
  \ar[d]<+6pt>_{t^{\underline{i}-i_{j_1}-...-i_{j_{k-m-2}}}_{i_{j_{k-m-1}}}}\\
                                C_{\underline{i}-i_{j_1}-...-i_{j_{k-m-1}}}
                                \ar[rrrr]^{j^{\underline{i}-i_{j_1}-...-i_{j_{k-m-1}}}_{i_{j_{k-m}}}}
                                \ar[drr]_{s^{\underline{i}-i_{j_1}-...-i_{j_{k-m-1}}}_{i_{j_{k-m}}}}&&&&
                                C_{\underline{i}-i_{j_1}-...-i_{j_{k-m-1}}}
                                \ar[dll]^{t^{\underline{i}-i_{j_1}-...-i_{j_{k-m-1}}}_{i_{j_{k-m}}}}\\
                                &&C_{\underline{i}-i_{j_1}-...-i_{j_{k-m}}} } $$
Also, as in \cite{Cam,cam-cubique} respectively for the globular geometry or for the
cubical geometry, it is possible to have a general notion of multiple $(\infty,m)$-structure :
this notion gives all possibilities of \textit{inverse structure} between the minimal $(\infty,m)$-structure and the maximal $(\infty,m)$-structure. Thus we define it as follows : a multiple set $(C,s,t)$ is equipped
with an $(\infty,m)$-structure if for all $k>m$, all $k$-color $\underline{i}$, and for all $1$-color $i_{j_1}\in S^1_{\underline{i}}$, there exist
an integer $q$ with $1\leq q\leq k-m-1$, there exist a $q$-color
$\{i_{j_1},...,i_{j_q}\}\in S^q_{\underline{i}}$, and there exist maps

$$\begin{tikzcd}
C_{\underline{i}}\arrow[rr,"j^{\underline{i}}_{i_{j_1}}"]&&C_{\underline{i}}
\end{tikzcd}, \qquad\begin{tikzcd}
C_{\underline{i}-i_{j_1}}\arrow[rr,"j^{\underline{i}-i_{j_1}}_{i_{j_2}}"]&&C_{\underline{i}-i_{j_1}}
\end{tikzcd}, \qquad\begin{tikzcd}
C_{\underline{i}-i_{j_1}-i_{j_2}}\arrow[rr,"j^{\underline{i}-i_{j_1}-i_{j_2}}_{i_{j_3}}"]&&
C_{\underline{i}-i_{j_1}-i_{j_2}}
\end{tikzcd},$$

$$\begin{tikzcd}
C_{\underline{i}-i_{j_1}-...-i_{j_{q-1}}}\arrow[rrr,"j^{\underline{i}-i_{j_1}-...-i_{j_{q-1}}}_{i_{j_q}}"]&&&
C_{\underline{i}-i_{j_1}-...-i_{j_{q-1}}}
\end{tikzcd}, \qquad\begin{tikzcd}
C_{\underline{i}-i_{j_1}-...-i_{j_{q}}}\arrow[rrr,"j^{\underline{i}-i_{j_1}-...-i_{j_{q}}}_{i_{j_{q+1}}}"]&&&
C_{\underline{i}-i_{j_1}-...-i_{j_{q}}}
\end{tikzcd}$$

such that we have the following diagrams in $\E$ which commute serially :

$$
  \xymatrixrowsep{8mm}
  \xymatrixcolsep{4mm}
  \xymatrix{
  C_{\underline{i}}\ar[rrrr]^{j^{\underline{i}}_{i_{j_1}}}
  \ar[d]<+1pt>_{s^{\underline{i}}_{i_{j_1}}}
  \ar[d]<+5pt>^{t^{\underline{i}}_{i_{j_1}}}&&&&C_{\underline{i}}
  \ar[d]<+3pt>^{s^{\underline{i}}_{i_{j_1}}}
  \ar[d]<+6pt>_{t^{\underline{i}}_{i_{j_1}}}\\
                                C_{\underline{i}-i_{j_1}}\ar[rrrr]^{j^{\underline{i}-i_{j_1}}_{i_{j_2}}}
    \ar[d]<+1pt>_{s^{\underline{i}-i_{j_1}}_{i_{j_2}}}
    \ar[d]<+5pt>^{t^{\underline{i}-i_{j_1}}_{i_{j_2}}}&&&&
      C_{\underline{i}-i_{j_1}}
      \ar[d]<+3pt>^{s^{\underline{i}-i_{j_1}}_{i_{j_2}}}
      \ar[d]<+6pt>_{t^{\underline{i}-i_{j_1}}_{i_{j_2}}}\\
                                C_{\underline{i}-i_{j_1}-i_{j_2}}
                                \ar[rrrr]^{j^{\underline{i}-i_{j_1}-i_{j_2}}_{i_{j_3}}}\ar@{.>}[d]<+1pt>^{}\ar@{.>}[d]<+5pt>^{}&&&&C_{\underline{i}-i_{j_1}-i_{j_2}}\ar@{.>}[d]<+3pt>^{}\ar@{.>}[d]<+6pt>_{}\\
                                C_{\underline{i}-i_{j_1}-...-i_{j_{q-1}}}
                                \ar[rrrr]^{j^{\underline{i}-i_{j_1}-...-i_{j_{q-1}}}_{i_{j_{q}}}}
  \ar[d]<+1pt>_{s^{\underline{i}-i_{j_1}-...-i_{j_{q-1}}}_{i_{j_{q}}}}
  \ar[d]<+5pt>^(.4){t^{\underline{i}-i_{j_1}-...-i_{j_{q-1}}}_{i_{j_{q}}}}&&&&
                                C_{\underline{i}-i_{j_1}-...-i_{j_{q-1}}}
  \ar[d]<+3pt>^{s^{\underline{i}-i_{j_1}-...-i_{j_{q-1}}}_{i_{j_{q}}}}
  \ar[d]<+6pt>_{t^{\underline{i}-i_{j_1}-...-i_{j_{q-1}}}_{i_{j_{q}}}}\\
                                C_{\underline{i}-i_{j_1}-...-i_{j_{q}}}
                                \ar[rrrr]^{j^{\underline{i}-i_{j_1}-...-i_{j_{q}}}_{i_{j_{q+1}}}}
                                \ar[drr]_{t^{\underline{i}-i_{j_1}-...-i_{j_{q}}}_{i_{j_{q+1}}}}&&&&
                                C_{\underline{i}-i_{j_1}-...-i_{j_{q}}}
                                \ar[dll]^{s^{\underline{i}-i_{j_1}-...-i_{j_{q}}}_{i_{j_{q+1}}}}\\
                                &&C_{\underline{i}-i_{j_1}-...-i_{j_{q+1}}} } $$

 $$
  \xymatrixrowsep{8mm}
  \xymatrixcolsep{4mm}
  \xymatrix{
  C_{\underline{i}}\ar[rrrr]^{j^{\underline{i}}_{i_{j_1}}}
  \ar[d]<+1pt>_{s^{\underline{i}}_{i_{j_1}}}
  \ar[d]<+5pt>^{t^{\underline{i}}_{i_{j_1}}}&&&&C_{\underline{i}}
  \ar[d]<+3pt>^{s^{\underline{i}}_{i_{j_1}}}
  \ar[d]<+6pt>_{t^{\underline{i}}_{i_{j_1}}}\\
                                C_{\underline{i}-i_{j_1}}\ar[rrrr]^{j^{\underline{i}-i_{j_1}}_{i_{j_2}}}
    \ar[d]<+1pt>_{s^{\underline{i}-i_{j_1}}_{i_{j_2}}}
    \ar[d]<+5pt>^{t^{\underline{i}-i_{j_1}}_{i_{j_2}}}&&&&
      C_{\underline{i}-i_{j_1}}
      \ar[d]<+3pt>^{s^{\underline{i}-i_{j_1}}_{i_{j_2}}}
      \ar[d]<+6pt>_{t^{\underline{i}-i_{j_1}}_{i_{j_2}}}\\
                                C_{\underline{i}-i_{j_1}-i_{j_2}}
                                \ar[rrrr]^{j^{\underline{i}-i_{j_1}-i_{j_2}}_{i_{j_3}}}\ar@{.>}[d]<+1pt>^{}\ar@{.>}[d]<+5pt>^{}&&&&C_{\underline{i}-i_{j_1}-i_{j_2}}\ar@{.>}[d]<+3pt>^{}\ar@{.>}[d]<+6pt>_{}\\
                                C_{\underline{i}-i_{j_1}-...-i_{j_{q-1}}}
                                \ar[rrrr]^{j^{\underline{i}-i_{j_1}-...-i_{j_{q-1}}}_{i_{j_{q}}}}
  \ar[d]<+1pt>_{s^{\underline{i}-i_{j_1}-...-i_{j_{q-1}}}_{i_{j_{q}}}}
  \ar[d]<+5pt>^(.4){t^{\underline{i}-i_{j_1}-...-i_{j_{q-1}}}_{i_{j_{q}}}}&&&&
                                C_{\underline{i}-i_{j_1}-...-i_{j_{q-1}}}
  \ar[d]<+3pt>^{s^{\underline{i}-i_{j_1}-...-i_{j_{q-1}}}_{i_{j_{q}}}}
  \ar[d]<+6pt>_{t^{\underline{i}-i_{j_1}-...-i_{j_{q-1}}}_{i_{j_{q}}}}\\
                                C_{\underline{i}-i_{j_1}-...-i_{j_{q}}}
                                \ar[rrrr]^{j^{\underline{i}-i_{j_1}-...-i_{j_{q}}}_{i_{j_{q+1}}}}
                                \ar[drr]_{s^{\underline{i}-i_{j_1}-...-i_{j_{q}}}_{i_{j_{q+1}}}}&&&&
                                C_{\underline{i}-i_{j_1}-...-i_{j_{q}}}
                                \ar[dll]^{t^{\underline{i}-i_{j_1}-...-i_{j_{q}}}_{i_{j_{q+1}}}}\\
                                &&C_{\underline{i}-i_{j_1}-...-i_{j_{q+1}}} } $$
and that the minimal $(\infty,m)$-structures are those with $q=1$, and the maximal $(\infty,m)$-structures are those with $q=k-m-1$ and where we have to consider all $q$-colors $\{i_{j_1},...,i_{j_q}\}\in S^q_{\underline{i}}$.
\end{remark}

\section{Weak multiple $(\infty,m)$-categories}
\label{weakm}

\subsection{Strict multiple $(\infty,m)$-categories}

A multiple $(\infty,m)$-magma $(M,s,t;j)$ is given by a multiple $\infty$-magma
$(M,s,t)$ equipped with an $(\infty,m)$-structure $j$ in the sense of \ref{reverses}, and
a morphism of $(\infty,m)$-magma is a morphism of $\infty\text{-}\mathbb{M}\mathbb{M}\text{ag}$
which is also a morphism of $(\infty,m)$-$\mathbb{M}\E$. The category of multiple $(\infty,m)$-magmas
is denoted $(\infty,m)\text{-}\mathbb{M}\mathbb{M}\text{ag}$. A reflexive multiple
$(\infty,m)$-magma $(M,s,t;1;j)$ is given by a multiple $\infty$-magma
$(M,s,t)$ equipped with reflexivity $1$ and equipped with an $(\infty,m)$-structure $j$. A
morphism of reflexive multiple $(\infty,m)$-magmas  is a morphism of
$\infty\text{-}\mathbb{M}_r\mathbb{M}\text{ag}$ which is also a morphism of $(\infty,m)$-$\mathbb{M}\E$. The category of reflexive multiple $(\infty,m)$-magmas
is denoted $(\infty,m)\text{-}\mathbb{M}_r\mathbb{M}\text{ag}$.

A strict multiple $(\infty,m)$-category is given by a strict multiple $\infty$-category
$(C,s,t,(\circ^{\underline{i}}_{\underline{i}-i_j})_{n\geq 1, \underline{i}\in\mathfrak{I}_n)}$
equipped with an $(\infty,m)$-structure. As for the globular geometry or the cubical geometry
(see \cite{Cam,cam-cubique}) it is not difficult to show that such $(\infty,m)$-structure is unique under this strictness. The category $\Mmc$ of strict multiple $(\infty,m)$-categories is the full subcategory
of $(\infty,m)\text{-}\mathbb{M}_r\mathbb{M}\text{ag}$ spanned by the strict multiple $(\infty,m)$-categories.

The third important monad of this article is the monad $\mathbb{S}^m=(S^m,\lambda^m,\nu^m)$ of strict multiple $(\infty,m)$-categories given by the forgetful functor

$$\begin{tikzcd}
\Mmc\arrow[rrr,"U"]&&&\mathbb{M}\E
\end{tikzcd}$$

and which is in fact monadic. The proof that $U$ is right adjoint and
monadic come from the underlying projective sketches of $\Mmc$ and $\mathbb{M}\E$ which are not difficult to be described, and it is enough to use the result in \ref{strict} the sketch of reversors described in \ref{reverses}, and then these results come from an easy application of the theorem of Foltz in \cite{foltz} and Lair in \cite{lair}.

\subsection{Multiple $(\infty,m)$-categorical stretchings}

A multiple $(\infty,m)$-categorical stretching
$(M,C,
(\pi_{\underline{i}})
_{n\geq 1,\underline{i}\in\mathfrak{I}_n},
([-,-]^{\underline{i}}_{\underline{i}+i_r})
_{n\geq 1,\underline{i}\in\mathfrak{I}_n})$
is given by the following datas :

\begin{itemize}

\item A reflexive multiple $(\infty,m)$-magma
$(M,s,t;1;j;(\circ^{\underline{i}}_{\underline{i}-i_j})_{n\geq 1, \underline{i}\in\mathfrak{I}_n)}$

\item A strict multiple $(\infty,m)$--category
$(C,s,t;1;j;(\circ^{\underline{i}}_{\underline{i}-i_j})_{n\geq 1, \underline{i}\in\mathfrak{I}_n)}$

\item A morphism of reflexive multiple $(\infty,m)$-magma
\begin{tikzcd}
M\arrow[rrrr,"\pi"]&&&&C
\end{tikzcd}

\item Operations $[-,-]^{\underline{i}}_{\underline{i}+i_r}$ :

$$\begin{tikzcd}
\underline{M}_{\underline{i}} \arrow[rrrr,"{[-;-]}^{\underline{i}}_{\underline{i}+i_r}"]&&&&M_{\underline{i}+i_r}
\end{tikzcd}$$
where $\underline{M}_{\underline{i}}=\{(\alpha,\beta)\in M_{\underline{i}}\times M_{\underline{i}} : \pi_{\underline{i}}(\alpha)=\pi_{\underline{i}}(\beta) \}$ such that (see \ref{stretchings}) :

\begin{itemize}

\item $s^{\underline{i}+i_r}_{i_s}[\alpha;\beta]^{\underline{i}}_{\underline{i}+i_r}=[s^{\underline{i}}_{i_s}(\alpha),s^{\underline{i}}_{i_s}(\beta)]^{\underline{i}-i_s}_{\underline{i}-i_s+i_r}$
and $t^{\underline{i}+i_r}_{i_s}[\alpha;\beta]^{\underline{i}}_{\underline{i}+i_r}=[t^{\underline{i}}_{i_s}(\alpha),t^{\underline{i}}_{i_s}(\beta)]^{\underline{i}-i_s}_{\underline{i}-i_s+i_r}$

\item $s^{\underline{i}+i_j}_{i_j}([\alpha;\beta]^{\underline{i}}_{\underline{i}+i_j})=\alpha$ and $t^{\underline{i}+i_j}_{i_j}([\alpha;\beta]^{\underline{i}}_{\underline{i}+i_j})=\beta$

\item $\pi_{\underline{i}+i_r}([\alpha;\beta]^{\underline{i}}_{\underline{i}+i_r})
=1^{\underline{i}}_{\underline{i}+i_r}(\pi_{\underline{i}}(\alpha))=1^{\underline{i}}_{\underline{i}+i_r}(\pi_{\underline{i}}(\beta))$

\end{itemize}

\end{itemize}

A morphism of multiple $(\infty,m)$-categorical stretchings
$$\begin{tikzcd}
\mathbb{E}\arrow[rrr,"{(m,c)}"]&&&\mathbb{E}'
\end{tikzcd}$$

is given by the following commutative square in
$(\infty,m)\text{-}\mathbb{M}_r\mathbb{M}\text{ag}$,

$$\begin{tikzcd}
M\arrow[rrr,"m"]
\arrow[dd,"\pi'"{left}]
&&&M'\arrow[dd,"\pi'"]\\\\
C \arrow[rrr,"c"{below}]&&&C'
\end{tikzcd}$$

thus we also have the following square in $\E$ for all $n$-color $\underline{i}\in\mathfrak{I}_n$
$(n\geq 1)$

$$\begin{tikzcd}
M_{\underline{i}}\arrow[rrr,"m_{\underline{i}}"]
\arrow[dd,"\pi_{\underline{i}}'"{left}]
&&&M'_{\underline{i}}\arrow[dd,"\pi'_{\underline{i}}"]\\\\
C_{\underline{i}} \arrow[rrr,"c_{\underline{i}}"{below}]&&&C'_{\underline{i}}
\end{tikzcd}$$

and we require for all $(\alpha,\beta)\in\underline{M}_{\underline{i}}$ the following
equality
$$m_{\underline{i}+i_r}([\alpha,\beta]^{\underline{i}}_{\underline{i}+i_r})=
[m_{\underline{i}}(\alpha),m_{\underline{i}}(\beta)]'^{\underline{i}}_{\underline{i}+i_r}$$

The category of multiple $(\infty,m)$-categorical stretchings is denoted $\MEtm$.

Now consider the forgetful functor :

$$\begin{tikzcd}
\MEtm\arrow[rrr,"U^m"]&&&\mathbb{M}\E
\end{tikzcd}$$

given by :

$$\begin{tikzcd}
(M,C,
(\pi_{\underline{i}})
_{n\geq 1,\underline{i}\in\mathfrak{I}_n},
([-,-]^{\underline{i}}_{\underline{i}+i_r})
_{n\geq 1,\underline{i}\in\mathfrak{I}_n})\arrow[rrr,mapsto]&&&M
\end{tikzcd}$$

\begin{proposition}
 The functor $U^m$ just above
 has a left adjoint which produces a monad
$\mathbb{M}^m=(M^m,\eta^m,\nu^m)$ on the category of multiple sets.
\end{proposition}

The proof that $U$ is right adjoint comes from the underlying projective sketches of $\MEtm$ and $\mathbb{M}\E$ which are not difficult to be described (see for example \cite{cam-cubique}). Then these results come from an application of the theorem of Foltz in \cite{foltz}.

\begin{definition}
A weak multiple $(\infty,m)$-category is a $\mathbb{M}^m$-algebra. The category of
weak multiple $(\infty,m)$-categories is denoted $\mathbb{M}^m\text{-}\mathbb{A}\text{lg}$.
Also models of weak multiple $\infty$-groupoids are given by the weak multiple
$(\infty,0)$-categories, and thus the category of weak multiple $\infty$-groupoids is denoted
$\mathbb{M}^0\text{-}\mathbb{A}\text{lg}$
\end{definition}

\bigbreak{}
\begin{minipage}{1.0\linewidth}
  Camell \textsc{Kachour}\\
  Institut de Recherche en Informatique Fondamentale\\
  CNRS (UMR 8243), Paris,
  France.\\
  Phone: 00 33645847325\\
  Email:\href{mailto:camell.kachour@gmail.com}{\url{camell.kachour@gmail.com}}
\end{minipage}

\end{document}